\documentclass[twocolumn, amsthm]{autart}
\usepackage{graphicx}
\usepackage[utf8]{inputenc}
\usepackage{mathptmx} % assumes new font selection scheme installed
\usepackage{newtxtext}
\usepackage{newtxmath}
\usepackage{textcomp}
\usepackage{enumitem}
\usepackage{dsfont, url}
\usepackage{mathtools}
\usepackage[usenames, dvipsnames]{xcolor}
\usepackage[%
colorlinks={true},%
citecolor={Blue},%
urlcolor={PineGreen},%
linkcolor={Sepia}%
]{hyperref}
\usepackage{float}
\usepackage{newfloat}
\usepackage[skip=1pt]{caption}
\usepackage{pgfplots}
\usepackage{amsmath} % assumes amsmath package installed
\usepackage{adjustbox}
\usepackage{verbatim} 
\usepackage{xargs}
\usepackage{cite}
\usepackage{todonotes}

\DeclareMathOperator*{\minimize}{minimize}
\DeclareMathOperator*{\argmin}{arg\ min}
\DeclareMathOperator*{\sbjto}{subject\ to}

\DeclareMathOperator{\rank}{rank}
\DeclareMathOperator{\trace}{tr}

\renewcommand{\leq}{\leqslant}

\renewcommand{\geq}{\geqslant}

\newcommand{\R}{\mathds{R}}
\newcommand{\Nz}{\mathds{N}_0}
\newcommand{\N}{\mathds{Z}_{+}}
\newcommand{\EE}{\mathds{E}}
\newcommand{\bmat}[1]{\begin{bmatrix}#1\end{bmatrix}}
\newcommand{\norm}[1]{\left|#1\right|}
\newcommand{\secref}[1]{\S \ref{#1}}

\renewcommand{\transp}{^\top}
\newcommand{\zeros}{\mathbf{0}}

\newcommand{\st}{x}

\newcommand{\stfilt}{\hat{x}}
\newcommand{\A}{A}

\newcommand{\meas}{y}
\newcommand{\C}{C}

\newcommand{\wnoise}{w}
\newcommand{\mnoise}{\varsigma}

\newcommand{\reachab}{\mathrm{R}}
\newcommand{\Let}{\coloneqq}
\newcommand{\teL}{\eqqcolon}

\theoremstyle{definition}
\newtheorem{theorem}{Theorem}
\newtheorem{lemma}{Lemma}

\newtheorem{remark}{Remark}
\newtheorem{definition}{Definition}
\newtheorem{example}{Experiment}
\newtheorem{assumption}{Assumption}

\allowdisplaybreaks

\begin{document}
	
	\begin{frontmatter}

		\title{Minimum variance constrained estimator} 
		
		\author[UIUC]{Prabhat K. Mishra} 
		\author[UIUC]{Girish Chowdhary}
		\author[UIUC]{Prashant G.\ Mehta}
		\address[UIUC]{Coordinated Science Laboratory,
			University of Illinois at Urbana Champaign (UIUC), USA.\tt\{pmishra, girishc, mehtapg\}@illinois.edu}
		

		\begin{abstract} 
				This paper is concerned with the problem of state estimation for
				discrete-time linear systems in the presence of additional (equality or
				inequality) constraints on the state (or estimate).  By use of the minimum variance
				duality, the estimation problem is converted into an optimal control
				problem.  Two algorithmic solutions are described: the full information
				estimator (FIE) and the moving horizon estimator (MHE).  The main
				result is to show that the proposed estimator is stable in the sense
				of an observer. The proposed algorithm is distinct from the standard
				algorithm for constrained state estimation based upon the use of the minimum energy duality.  The two
				are compared numerically on the benchmark batch reactor process model.
		\end{abstract}
		
	\end{frontmatter}

	\section{Introduction}
	\label{sec:introduction}
	\par In many practical estimation problems arising in control applications,
	there are invariably additional constraints on the state process~\cite{constrained_control_estimation}.  In such applications,
	Kalman filter (KF) may yield sub-optimal estimates that violate the
	constraints. It is notable also that the KF  is
	derived under the assumption of (unbounded) Gaussian noise, which is
	also unrealistic in the constrained settings of the problem. In
	particular, in the presence of unbounded noise, local stability
	results are not applicable and global stability results are very
	conservative due to actuator saturation \cite{ref:ChaRamHokLyg-12,
		PDQ-LCSS}.  Although clever modifications in KF are still
	possible~\cite{KF_constraints}, 
	the stability and optimality properties of such modifications require
	further investigation~\cite{KF_constraints_survey}.  For these
	reasons, constrained estimation is a problem of paramount practical
	importance; c.f.,~\cite{constrained_control_estimation} for a book
	length treatment.	
	\par A popular strategy for constrained estimation is based on the use of
	duality between estimation and optimal control.  A practical advantage
	of converting a constrained estimation problem into a constrained
	optimal control problem is that model predictive control (MPC)
	methods, algorithms, and softwares can readily be applied to obtain a
	solution.  The resulting estimation algorithms are referred to as the
	\emph{full information estimator} (FIE), when {\em all} the observations are 
	used, and is \emph{moving horizon estimator} (MHE), when a moving window of most
	recent past observations are used.  Practically, a MHE algorithm is
	preferred because the number of decision variables in the optimization
	problem do not increase as more observations are collected.         
	
	\par In linear settings of the problem, there are in general two types of duality: the
	minimum energy (or maximum likelihood) duality and the minimum variance duality. Refer to \cite{pearson66_duality, duality_variational, duality_general} for more discussion on duality. For the
	construction of estimators, the minimum energy duality is by far the
	more popular technique with contributions
	in~\cite{Copp_Hespanha_simultaneous,coupled_methods,
		farina_distributedMHE, mhe15_interacting, MHE_Switched,
		MHE_Brembeck, fast_MHE} and numerical algorithms in~\cite{mhe_numerical, Findeisen_MHEcomputation}.  Although minimum
	variance control has attracted much
	attention~\cite{bakolas_min_variance, bakolas18_covarianceControl,
		covariance_steering_games}, and these recent papers provide
	motivation also for our work, the use of minimum variance duality for constrained estimation
	has received comparatively less attention.   
	
	\par State estimation problem for linear systems with equality constraints
	is considered in \cite{equality_constrainted_estimator,
		equality_constrainted_Bernstein}, and with inequality constraints in
	\cite{inequality_LS}. Since the number of decision variables in the
	underlying optimization program increases as more
	measurements are collected, an MHE algorithm is proposed
	in the early work of~\cite{MHE68} in the absence of constraints.  
	This algorithm is extended in
	\cite{MHE93} to incorporate constraints. The stability properties of
	this constrained MHE algorithm are studied rigorously in
	\cite{MHE01}.  Enhancements of these basic algorithms have been
	considered both in deterministic (observer design)~\cite{MHE03,MHE10,
		MHE2014_constrained, MHE2018_proximity} and stochastic (filter
	design)~\cite{MHE68, MHE93, MHE01, MHE02} settings of the problem.
	In stochastic settings, the MHE optimal control problem is still
	deterministic but statistical information about uncertainties and
	prior are used to design the (maximum likelihood-type) objective
	function.  More recent extensions include a game theoretic formulation
	in~\cite{MHE18_game}. It is noted that \cite{MHE68, MHE93, MHE01, MHE02, MHE18_game} are
	based on the use of the minimum energy duality. We refer readers to \cite[Appendix B]{rao_thesis} for a quick
	review of duality and to \cite{mle68} for minimum energy
	duality in particular.   	
	\par In this paper, an alternate form of duality, viz., the minimum
	variance duality is employed to transform the minimum variance
	estimation problem into a deterministic optimal control problem. The
	state estimate is constructed as a linear function of past measurements.
	Without constraints, the optimal estimate is equivalent to a Kalman
	filter. Both the FIE and the MHE are described for the unconstrained
	case, together with expression for choosing the terminal cost in the
	MHE.  
	\par The main focus of this paper is on the modification of these
	(unconstrained) FIE and MHE algorithms in the presence of constraints.
	In particular, a certain approximate expression for the terminal cost
	is introduced for the constrained MHE.  The main result of this paper
	is to establish sufficient conditions to obtain stability (in
	the sense of an observer) for the constrained FIE and MHE algorithms.
	Furthermore, we also establish a certain type of stochastic stability
	by showing that the variance of the constrained FIE converges under certain
	technical conditions.
	\par Although estimators based on minumum variance duality are
	less well studied \cite{kim2019_CDC}, some closely related estimators have appeared in ~\cite{RHKF99, mv_unbiased_inpput, MHE02,mvFIR07, min_variance_fir_tv, mv_ui}. In
	contrast to our paper, these prior works do not incorporate equality
	or inequality constraint on state (or estimate) in the estimator
	design. The original contributions of our paper are as follows:
\begin{itemize}[leftmargin = *]
		\item Based on minimum variance duality, a MHE is presented in
		\eqref{e:MHE} and its equivalence with FIE is shown in Lemma
		\ref{lem:uncsontrained_mhe}. This contribution is different from
		\cite{MHE02, mvFIR07} in the sense that unbiasedness constraints are
		not required. The proposed estimator \eqref{e:estimator_mhe} is equivalent to KF.
		\item Constrained FIE and MHE algorithms are presented in
		\eqref{e:constrained_fie} and \eqref{e:constrained_MHE},
		respectively. Apart from the fact that these algorithms are distinct
		from~\cite{MHE01}, our minimum variance-based approach has certain
		technical advantages.  
		\item Although the notion of stability is borrowed from~\cite{MHE01},
		Theorems \ref{prop:FIE_stability} and \ref{prop:MHE_stability} are
		first such results on stability of constrained minimum variance
		estimators.  
		\item Under certain technical conditions, the variance of constrained
		FIE is shown to converge in Theorem
		\ref{th:stochastic_cfi}. 
\end{itemize}
	\par The remainder of this paper is organized as follows: The problem
	statement appears in \secref{s:problem statement} followed by a
	description of the minimum variance duality for the construction of
	the unconstrained estimators, both FIE and MHE, in
	\secref{s:unconstrained}.  These are extended to the constrained case
	in \secref{s:constrained}.  The main results on stability of the
	constrained FIE and MHE appear in \secref{s:stability}.
	The algorithms are illustrated with the aid of some numerical
	experiments in \secref{s:experiments}. The paper closes with some
	conclusions and directions for future research in \secref{s:conclusion}. All the proofs appear as part
	of the two appendices,\secref{s:proofs_unconstrained} and
	\secref{s:proofs_constrained}, for the unconstrained and the
	constrained cases, respectively.  
	\par Let $\R, \Nz, \N$ denote the set of real numbers, the non-negative integers and the positive integers, respectively. We use the symbols \( \zeros \) and \(I\) to denote zero matrix and identity matrix, respectively, of appropriate dimensions.  
	For any vector or matrix sequence \((M_n)_{n\in\Nz} \in \R^{r\times m}, r,m \in \N\), let \( M_{n:k} \in \R^{rk\times m} \) denote the matrix \(\bmat{M_n\transp & M_{n+1}\transp & \cdots & M_{n+k-1}\transp}\transp\), \(k\in\N\). 
	Let $\lambda_{\max}(M)$ denote the largest eigenvalue value of $M$, $\lambda_{\min}(M)$ its smallest eigenvalue, $M^\dagger$ its Moore-Penrose pseudo inverse and $\trace(M)$ its trace. 
	The Euclidean norm of a vector $A$ is denoted by $\norm{A}$. The Frobenius norm of a matrix $A$ is denoted by $\norm{A}_F$. A $t$ step reachability matrix of a matrix pair $(A,B)$ is given by $\reachab_t(A,B) \Let \bmat{A^{t-1}B & \ldots & AB & B}$.
	%%%%%%%%%%%%%%%%%%%%%%%%%%%%%%%%%%%%%%%%%%%%%%%%%%%%%%%%%%%%%%%%%%%%%%%%%%%%%%%%%%%%%%%%%%%
	\section{Problem statement}\label{s:problem statement}
	Consider a linear discrete-time system
	\begin{equation}\label{e:process}
	\begin{aligned}
	\st_{t+1} &= \A \st_t + \wnoise_t,\\
	\meas_t &= \C \st_t + \mnoise_t,
	\end{aligned}
	\end{equation}	 
	where $\st_t \in \bar{\mathcal{X}} \subset \R^d$, $\meas_t \in \R^q$ are state and measurement of the system at time $t$, respectively. The system matrix $A \neq \zeros$. The additive process noise $\wnoise_t$ and the measurement noise $\mnoise_t$ are mean zero, mutually independent and identically distributed random vectors with variance $Q$ and $R$, respectively.  The initial state of the system $\st_0$ is a random vector with mean $\stfilt_0^-$ and variance $\Sigma_0^-$, and is independent of the process noise and the measurement noise. 
	
	\par The minimum variance estimation problem is to compute
	$\stfilt_t$ at time $t$ such that the variance of error $\st_t - \stfilt_t$ is minimized over some class of admissible estimators. In
		this paper, the admissible estimators are assumed to be linear
		deterministic functions of available measurements.  
		It is also assumed that some additional insight into the states (or estimates) is given in terms of equality and inequality constraints such that the estimated states belong to a convex set $\mathcal{X} \supseteq \bar{\mathcal{X}}$ , i.\ e.\ , $\stfilt_t \in \mathcal{X}$ for all $t$. We make the following assumption:
		\begin{assumption}\label{as:positive_invariance}
			The set $\mathcal{X}$ is positively invariant under the nominal dynamics, i.\ e.\ , 
			$\A x \in \mathcal{X} \text{ for every } x \in \mathcal{X}$. 	
		\end{assumption} 
		The above asssumption is meaningful.
		Suppose $A\stfilt_t \notin \mathcal{X}$ for some $\stfilt_t \in
		\mathcal{X}$ then there is a non-zero probability that $\st_{t+1}
		\notin \bar{\mathcal{X}}$ for random $\wnoise_t$, e.g., when
		$\stfilt_t = \st_t$ and a bounded disturbance set with known bounds is
		not safely prescribed. The optimization problem is as follows:
	\begin{equation}\label{e:estimation_problem}
	\min_{\stfilt_t \in \mathcal{X}}\quad  \EE \left[ \norm{\st_t - \stfilt_t}^2 \right].
	\end{equation} 
	The solution approach is based on duality between estimation
	and control. In the following section, we begin by presenting
	an unconstrained estimator which is useful for the development
	of a constrained estimator in
	\secref{s:constrained}.
	\section{Minimum variance estimators}
	\subsection{Unconstrained estimator}\label{s:unconstrained}
	In this section, we assume $\mathcal{X} = \R^d$, i.e., the constraints are not present. We are interested in an estimator linearly parameterized in the innovation terms as follows:
		\begin{equation}\label{e:stfilt without z}
		\stfilt_t = A^t\stfilt_0^- -\sum_{i=0}^t \alpha_i \transp (\meas_{t-i} - C A^{t-i}\stfilt_0^-),
		\end{equation}
		in which weights $\alpha_i \in \R^{q\times d}$ are the
		decision variables for the optimization problem~\eqref{e:estimation_problem}.
		In order to convert the minimum variance estimation objective
		into an optimal control problem, a dual process (in forward
		time) is introduced:
		\begin{equation}\label{e:dual_process}
		\begin{aligned}
		z_{i+1} &= \A\transp z_i + \C \transp \alpha_{i+1}; \quad i = 0, \ldots, t-1, \\
		z_0 &= I + C\transp \alpha_0, 
		\end{aligned}
		\end{equation}
		where $z_i \in \R^{d\times d}$ is a matrix valued dual state and $\alpha_i \in \R^{q \times d}$ is control signal for the dual process. From \eqref{e:dual_process} we have
		\begin{equation}\label{e:dual_state}
		z_t\transp = A^t + \sum_{i=0}^t \alpha_{i}\transp C A^{t-i}.
		\end{equation}
		By substituting \eqref{e:dual_state} into \eqref{e:stfilt without z}, we get the following expression:
		\begin{equation}\label{e:estimator}
		\stfilt_t = z_t\transp \stfilt_0^- - \sum_{i=0}^{t} \alpha_i \transp \meas_{t-i}.
		\end{equation}
A slight modification of the standard result on minimum variance duality \cite[Page 238]{Astrom_stochastic_control} \footnote{See \cite[Exercise 1, Page 240]{Astrom_stochastic_control} in which invertibility of the system matrix $A$ is assumed to define a dual process.}, in which only the past measurements are used to design an estimator, i.\ e.\ $\alpha_t = \zeros$, is required to include the current measurement. Let $\ell_i \Let z_i\transp Q z_i + \alpha_i\transp R\alpha_i$, $\Gamma_0(z_t) \Let z_t \transp \Sigma_{0}^- z_t$ and 	
	\begin{equation}\label{e: matrix without terminal cost}
	S_t(\alpha_{0:t+1}) \Let \alpha_t\transp R\alpha_t + \sum_{i=0}^{t-1} \ell_i. 
	\end{equation}
	The estimate \eqref{e:stfilt without z} takes into account all
	measurements available at time $t$. Therefore, the
	corresponding estimator is called \emph{full information
		estimator} (FIE). Using the dual process
	\eqref{e:dual_process}, the FIE optimal control problem is
	expressed as follows:
	\begin{equation}\label{e:fie}
	\text{FIE:} \quad \left\{
	\begin{aligned}
	\minimize_{\alpha_{0:t+1}} & \quad \trace(\Gamma_0(z_t) + S_t(\alpha_{0:t+1})) \\
	\sbjto & \quad \text{dual dynamics } \eqref{e:dual_process}.	 
	\end{aligned}
	\right.	
	\end{equation}
	FIE \eqref{e:fie} is solved at each time $t =0, 1,
	\ldots $. The resulting optimal solution is denoted as
	$\alpha_{0:t+1\mid t}$, where $\alpha_{k\mid t}$ is the
	optimal weight $\alpha_k$ computed at time $t$. Set
	\begin{equation}\label{e:variance_unconstrained}
	\Sigma_t \Let \Gamma_0(z_{t\mid t}) + S_t(\alpha_{0:t+1\mid t}),
	\end{equation}
	where $S_t(\alpha_{0:t+1\mid t})$ is the optimal value of
	$S_t(\alpha_{0:t+1})$ obtained by solving FIE
	\eqref{e:fie}. Then the optimal value of the
	objective function in \eqref{e:fie} is
	$\trace(\Sigma_t)$. The estimate $\stfilt_{t\mid
		t}$ at time $t$ is obtained by substituting the optimal
	values $\alpha_{0:t+1\mid t}$ and $z_{t\mid t}$ in
	\eqref{e:estimator}. In the remainder of the manuscript, we
	will use $\stfilt_t$ to denote the estimate obtained by
	substituting the optimizers in \eqref{e:estimator}. 
	We have the following Lemma to show the equivalence of FIE \eqref{e:fie} and \eqref{e:estimation_problem} whenever $\mathcal{X}=	\R^d$.  	
	\begin{lemma}\label{lem:duality}
		Consider the system \eqref{e:process} and the dual process \eqref{e:dual_process}. If $\stfilt_t$ is given by \eqref{e:estimator} then 
		\[ \EE \left[ \norm{\st_t - \stfilt_t}^2 \right] = \trace(\Gamma_0(z_t) + S_t(\alpha_{0:t+1})). \]
	\end{lemma}

	\begin{remark}
		\rm{
			The dual process is typically considered
			backward in time.  However, because the
			optimal control problem is deterministic, a
			forward time dual process may equivalently be considered
			simply by renaming the indices.  This is done
			here to yield the standard form of an optimal
			control problem where the time arrow is forward.
		}
	\end{remark}
	\par We present a finite horizon approximation
	of FIE \eqref{e:fie}, which we refer to as \emph{moving horizon
		estimator} (MHE). For this purpose, define
	\begin{equation}\label{e:prior}
	\Sigma_t^- \Let A \Sigma_{t-1} A \transp + Q \text{ and } \stfilt_t^- \Let A\stfilt_{t-1}.
	\end{equation}
	Fix $N \in \Nz$ and for $t \geq N+1$ define 
	\begin{equation}\label{e:terminal_cost_unconstrained}
	\Gamma_{t-N}(z_N) \Let z_N\transp \Sigma_{t-N}^- z_N .
	\end{equation}     
	The unconstrained MHE is
	as follows:
	\begin{equation} \label{e:MHE}
	\text{MHE:} \quad \left\{
	\begin{aligned}
	\minimize_{\alpha_{0:N+1}} & \quad \trace(\Gamma_{t-N}(z_N) + S_N(\alpha_{0:N+1})) \\
	\sbjto & \quad \text{dual dynamics } \eqref{e:dual_process}.	 	
	\end{aligned} 
	\right.
	\end{equation}
	For $t\leq N$, set $\Gamma_{t-N} = \Gamma_0, S_N = S_t$, which
	is identical to solving the FIE problem \eqref{e:fie}. For
	$t\geq N+1$, the MHE problem utilizes the most recent $N+1$
	measurements together with the previously computed
	$\Sigma_{t-N-1}$ to obtain $\Sigma_{t-N}^-$. The resulting estimator and the error covariance matrix are 
	\begin{align}
	\stfilt_t &= z_{N\mid t}\transp \stfilt_{t-N}^- - \sum_{i=0}^{N}\alpha_{i\mid t}\transp \meas_{t-i}, \label{e:estimator_mhe}\\
	\Sigma_t &= \trace(\Gamma_{t-N}(z_{N\mid t}) + S_N(\alpha_{0:N+1 \mid t})), \label{e:variance_mhe}
	\end{align}
	where $\alpha_{i\mid t}$ for $i=0,\ldots ,N$, and $z_{N\mid
		t}$ are obtained by solving MHE \eqref{e:MHE} at time $t$.
	It is straightforward to show that, when $N=0$, MHE
	\eqref{e:MHE} is the KF.  A direct implication of
	dynamic programming is the following result:
	\begin{lemma}\label{lem:uncsontrained_mhe} 
		If $R \succ 0$ then FIE \eqref{e:fie} is equivalent to MHE \eqref{e:MHE} and the estimate \eqref{e:estimator} is equal to the estimate \eqref{e:estimator_mhe}.
	\end{lemma}	
Proofs of lemmas \ref{lem:duality} and \ref{lem:uncsontrained_mhe} are given in the Appendix \ref{s:proofs_unconstrained}. 

	\subsection{Constrained estimator}\label{s:constrained}
	If the matrix pair $(A,C)$ is observable then there exists an
	integer $n \leq d \in \N$ such that
	$\rank(\reachab_n(A\transp, C\transp)) = d$.  The smallest
	such $n$ is referred to as the observability index of
	$(A,C)$. Our construction of the constrained FIE depends on
	$n$. In particular, we augment the FIE~\eqref{e:fie} with the
	following additional constraints:
	\begin{equation}\label{e:intermediate constraint}
	\begin{aligned}
	& z_{t-j}\transp \stfilt_0^- - \sum_{i=0}^{t-j}\alpha_i\transp \meas_{t-j-i} \in \mathcal{X},
	\end{aligned}
	\end{equation}
	where  $j=0$ for $t \leq n$ and $j =  0, \ldots , t-n$, for $t \geq n+1$. Note that the left hand side of the
	constraint is same as $\stfilt_{t-j}$ according
	to~\eqref{e:estimator}. 
	Although we are interested in this constraint only with $j=0$,
	inclusion of the intermediate constraints, for $j=1, \ldots
	,t-n$, helps to ensure some properties. Additional details on this
	appear in the next section. The constrained FIE problem is
	formally defined as follows:
	\begin{equation}\label{e:constrained_fie}
	\text{CFIE:} \quad \left\{
	\begin{aligned}
	\minimize_{\alpha_{0:t+1}} & \quad \trace(\Gamma_0(z_t) + S_t(\alpha_{0:t+1})) \\
	\sbjto & \quad \text{dual dynamics } \eqref{e:dual_process},	 \\
	& \quad \text{constraints } \eqref{e:intermediate constraint}.
	\end{aligned}
	\right.	
	\end{equation}	
	The solution of the CFIE \eqref{e:constrained_fie} is
	used to construct the constrained full information estimate by
	using the right hand side of \eqref{e:estimator}.  It is denoted $\stfilt_t^{\text{cf}}$ to distinguish it from unconstrained estimate $\stfilt_t$ obtained by solving \eqref{e:fie} or \eqref{e:MHE}. In particular,
	\begin{equation}\label{e:cfie_mean_variance}
	\begin{aligned}
	\stfilt_t^{\text{cf}} \Let z_{t\mid t}\transp \stfilt_0^- - \sum_{i=0}^{t} \alpha_{i\mid t} \transp \meas_{t-i},\\
	\Sigma_t^{\text{cf}} \Let  \Gamma_0(z_{t\mid t}) + S_t(\alpha_{0:t+1\mid t}),
	\end{aligned}
	\end{equation}
	where $z_{t\mid t}$ and $\alpha_{0:t+1\mid t}$ are obtained by solving \eqref{e:constrained_fie}. 	
		\begin{remark}[Feasibility and convexity]\label{rem:feasibility}
			If $\stfilt_0^{-} \in \mathcal{X}$ then the optimal control problem \eqref{e:constrained_fie} is feasible for all $t$ because $\alpha_{0:t+1} = \zeros$ satisfies \eqref{e:intermediate constraint}. The left hand side of \eqref{e:intermediate constraint} is affine in decision variables $\alpha_{0:t+1}$ and the set $\mathcal{X}$ is convex. The set of decision variables $\alpha_{0:t+1}$ in \eqref{e:intermediate constraint} is convex due to the fact that the inverse image of a convex set under an affine function is convex \cite[Page 38]{boyd_convex_optimization} and the intersection of convex sets is convex.
		\end{remark}  	
	\begin{remark}
		\rm{
			The right hand side of \eqref{e:estimator} is
			linear in the past measurements.  The justification comes from the 
			unconstrained linear Gaussian case where such
			a structure is sufficient to obtain the
			minimum variance estimator.  In the presence
			of constraints and non-Gaussian noise, an
			optimal estimate may not be linear in the past
			measurements.  It is noted that the assumed
			structure is also nonlinear because of the dependance of $\alpha_{0:t+1}$ on $\meas_{0:t+1}$ via constraint~\eqref{e:intermediate constraint}. 	
		}
	\end{remark}
	In the presence of constraints, the design of an MHE algorithm, that
	is provably equivalent to the FIE algorithm, is challenging because
	of the difficulty in approximating the terminal cost.  Therefore,
	approximation of the terminal cost (which is also referred to as
	\emph{arrival cost} in the standard MHE literature) is necessary.  The
	goal is to approximate the FIE as closely as possible while
	maintaining computational tractability and guaranteeing stability.
	
	\par Similar to CFIE \eqref{e:constrained_fie}, constrained MHE can also be defined by
	adding extra constraints to the unconstrained MHE
	\eqref{e:MHE}.  
	The constrained MHE estimator is denoted as
	$\stfilt_t^{\text{cm}}$, where the superscript cm is used to
	reflect the fact that this estimate at time $t$ may be different from the unconstrained estimate $\stfilt_t$ and the CFIE estimate $\stfilt_t^{\text{cf}}$. Similarly, the corresponding error covariance matrix is denoted by $\Sigma_t^{\text{cm}}$ to distinguish it from \eqref{e:variance_mhe}. 
	\par We need to define priors $\Sigma_{t-N}^{\text{cm}-}$ and $\stfilt_{t-N}^{\text{cm}-}$ to compute the terminal cost of the constrained MHE and its estimate as we did in \eqref{e:terminal_cost_unconstrained} and \eqref{e:estimator_mhe}, respectively, for the unconstrained case. One possible choice is to use $\Sigma_{t-N}^-$ and $\stfilt_{t-N}^-$ obtained from the unconstrained case by using \eqref{e:prior} and \eqref{e:MHE}, which is same as running a KF in parallel. The standard MHE \cite{MHE01} follows this approach. Other MHE approaches like \cite{MHE2014_constrained} also use priors from the unconstrained case. Since our approach not only gives an estimated state which satisfies constraints but also an error covariance matrix, it is reasonable to replace $\Sigma_{t-1}$ in \eqref{e:prior} by $\Sigma_{t-1}^{\text{cm}}$ to get $\Sigma_{t}^{\text{cm}-}$ and $\stfilt_{t-1}$ by $\stfilt_{t-1}^{\text{cm}}$ to get $\stfilt_t^{\text{cm}-}$. This choice is intuitive because the pair $(\stfilt_t^{\text{cm}-}, \Sigma_{t}^{\text{cm}-})$ represents our prior knowledge about the pair $(\stfilt_t^{\text{cm}}, \Sigma_{t}^{\text{cm}})$ in the presence of constraints. More precisely,
	\begin{equation}\label{e:prior constrained}
	\Sigma_t^{\text{cm}-} \Let A \Sigma_{t-1}^{\text{cm}} A \transp + Q \text{ and } \stfilt_t^{\text{cm}-} \Let A\stfilt_{t-1}^{\text{cm}}.
	\end{equation}
	\par The constrained MHE problem is formally written as follows:
	\begin{equation} \label{e:constrained_MHE}
	\text{CMHE:} \quad \left\{
	\begin{aligned}
	\minimize_{\alpha_{0:N+1}} & \quad \trace(\Gamma_{t-N}^{\text{cm}}(z_N) + S_N(\alpha_{0:N+1})) \\
	\sbjto & \quad \text{dual dynamics } \eqref{e:dual_process}, \text{ and } \\
	& \quad  z_N\transp\stfilt_{t-N}^{\text{cm}-} - \sum_{i=0}^{N} \alpha_i \transp \meas_{t-i}  \in \mathcal{X},  	
	\end{aligned} 
	\right.
	\end{equation}
	where $\Gamma_{t-N}^{\text{cm}}(z_N) \Let z_N \transp \Sigma_{t-N}^{\text{cm}-} z_N $ for $t\geq N+1$. Similar to MHE \eqref{e:MHE} for $t\leq N$, we set $\Gamma_{t-N}^{\text{cm}} = \Gamma_0, S_N = S_t$ and similarly modify constraint by taking all $t+1$ measurements. Alternatively, we can run CFIE \eqref{e:constrained_fie} for $t \leq N$.
	Further, the estimate \eqref{e:estimator} and corresponding covariance matrix can be written as	
	\begin{equation}\label{e:constrained_estimator}
	\begin{aligned}
	\stfilt_t^{\text{cm}} \Let z_{N\mid t}\transp \stfilt_{t-N}^{\text{cm}-} - \sum_{i=0}^{N}\alpha_{i\mid t}\transp \meas_{t-i},\\
	\Sigma_t^{\text{cm}} \Let  \Gamma_{t-N}^{\text{cm}}(z_{N\mid t}) + S_N(\alpha_{0:N+1\mid t}),
	\end{aligned}
	\end{equation}
	where $\alpha_{i\mid t}$ for $i=0,\ldots ,N$, and $z_{N\mid t}$, are obtained by solving CMHE problem \eqref{e:constrained_MHE} at time $t$. 
	\section{Main results}\label{s:stability}
	\par In this section, stability of the proposed constrained estimators is presented by using the notion of stability introduced in \cite{MHE01}. Recall that the classical notion of stability of an observer is obtained by modifying the definition of the stability of a regulator. In an analogous manner, the definition of the stability
	of a constrained regulator, which is given in \cite[\S 2]{Keerthi-Gilbert88}, is modified in \cite{MHE01} to introduce the following definition:
	\begin{definition}[\cite{MHE01, Keerthi-Gilbert88}]\label{def:stable_observer}
		The estimator is a stable observer for the system
		\begin{align}\label{e:nominal_system}
		\st_{t+1} &= A\st_t; \quad \meas_t = C\st_t; \quad \st_t \in \mathcal{X},
		\end{align}
		if for any $\varepsilon > 0$, there exists $\delta > 0$ and $T \in \N$ such that if $\stfilt_0^- \in \mathcal{X}$ and $\norm{\st_0 - \stfilt_0^-} \leq \delta$ then $\norm{\stfilt_t - A^t\st_0} \leq \varepsilon$ for all $t \geq T$. If in addition, $\stfilt_t \rightarrow A^t\st_0$ as $t \rightarrow \infty$ then the estimator is called asymptotically stable observer for the system \eqref{e:nominal_system}.	
	\end{definition}	
	Our approach has a minor advantage over \cite{MHE01} in the
	sense that a key assumption is relaxed. In particular, we do
	not assume any upper bound on cost a priori but it comes
	naturally from the observability of the system. For the
	stability of CFIE we need one of the following two conditions
	to hold:
	\begin{enumerate}[label={(C\arabic*}), leftmargin=*, widest=3, align=left, start=1]
		\item \label{as:larger_Q} $Q - \Sigma_{0}^- \succeq 0 $.
		\item \label{as:stabilizing_condition} There exists some $K_t \in \R^{q\times d}$ at each time $t \geq n+1$ such that $\alpha_t = K_t z_{t-1 \mid t-1}$ satisfies \eqref{e:intermediate constraint} for $j=0$ and the following stability criterion with $\tilde{A}_t = A\transp + C\transp K_t$:
		\begin{equation}\label{e:monotonicity_Lyapunov_condition}
		\tilde{A}_t\transp \Sigma_0^- \tilde{A}_t - \Sigma_0^- \preceq - (K_t\transp R K_t + Q).
		\end{equation}
	\end{enumerate}	
	The main stability result for the CFIE is as follows:
	\begin{theorem}\label{prop:FIE_stability}
		Suppose Assumption \ref{as:positive_invariance}
		holds, $\norm{\stfilt_0^- - \st_0} < \infty$, $\Sigma_0^- \succ 0$, $(A,C)$ is observable,
		and one of the two conditions, either
		\ref{as:larger_Q} or \ref{as:stabilizing_condition},
		is satisfied.  Then CFIE is an asymptotically stable observer for the system \eqref{e:nominal_system}. 
	\end{theorem}
	\begin{remark}
		It is easily verified that the conditions \ref{as:larger_Q}
		and \ref{as:stabilizing_condition} can not simultaneously hold
		unless $A = \zeros$, which because the matrix pair $(A,C)$ is
		observable, represents a trivially false case. Let, if possible, \ref{as:larger_Q} and \ref{as:stabilizing_condition} hold simultaneously then \ref{as:stabilizing_condition} gives
		\begin{equation}\label{e:stabilizing_condition_implication}
		0 \preceq \tilde{A}_t\transp \Sigma_0^- \tilde{A}_t  \preceq - (K_t\transp R K_t + Q - \Sigma_0^-) \preceq 0,
		\end{equation} 
		which implies $\tilde{A}_t\transp \Sigma_0^- \tilde{A}_t =
		K_t\transp R K_t + Q - \Sigma_0^- = \zeros$. Therefore,
		$\tilde{A}_t = \zeros$ because $\Sigma_0^- \succ 0$ and
		$K_t\transp R K_t + Q  = \Sigma_0^-$, which results in $Q
		\preceq \Sigma_0^-$ and due to \ref{as:larger_Q} we get $Q =
		\Sigma_0^-$. By substituting $Q=\Sigma_0^-$ in
		\eqref{e:stabilizing_condition_implication}, we get
		$K_t\transp R K_t = \zeros$, which results in $K_t = \zeros$
		because $R \succ 0$. Since $\tilde{A}_t = \zeros$ due to
		\eqref{e:stabilizing_condition_implication}, the substitution
		of $K_t = \zeros$ shows that $A = \zeros$.   
	\end{remark}
	We have the following result on stability of CMHE:
		\begin{theorem}\label{prop:MHE_stability}
			Suppose Assumption \ref{as:positive_invariance}
			holds, $\Sigma_0^- \succ 0$, $R \succ 0$ and $(A,C)$ is observable then for $N\geq n$, CMHE is stable observer for the system \eqref{e:nominal_system}. If, in addition, $Q \succ 0$, $\norm{\stfilt_0^- - \st_0} < \infty$, then CMHE is asymptotically stable observer for the system \eqref{e:nominal_system}.  
		\end{theorem}
		In theorems \ref{prop:FIE_stability} and \ref{prop:MHE_stability}, we proved stability of the proposed estimators in the sense of an observer. Since the cost function represents variance in the proposed approach, we get its convergence for the system \eqref{e:process} also under the following assumption:
		\begin{assumption}\label{as:sufficient_feasibility}
			There exist $\alpha_0 \in \R^{q \times d}$, and a sequence of matrices $(K_i)_{i\in \N}$ such that $\alpha_{i+1} = K_i z_i$ and $\alpha_0$ satisfy \eqref{e:intermediate constraint}. There exist $\lambda_0 > 0$, $\lambda_i < 1$ for $i \in \N$ such that
			\begin{enumerate}[label={(C\arabic*}), leftmargin=*, widest=3, align=left, start=3]
				\item \label{as:upperbound_1} $(I+C\transp \alpha_0)\transp Q (I+C\transp \alpha_0)+ \alpha_0\transp R \alpha_0 \leq \lambda_0 Q$
				\item \label{as:upperbound_2} $(A\transp + C\transp K_i)\transp Q (A\transp + C\transp K_i) + K_i \transp R K_i \leq \lambda_i Q \text{ for } i \in \N$
			\end{enumerate}	
		\end{assumption}
		The above assumption gives a sufficient condition for the feasibility of \eqref{e:constrained_fie} and the existence of a stabilizing controller for the dual process \eqref{e:dual_process}. Notice that \eqref{e:constrained_fie} is feasible due to the Remark \ref{rem:feasibility}. The above assumption helps us to get an upper bound of the cost in \eqref{e:constrained_fie}.  	
		We have the following result: 
		\begin{theorem}\label{th:stochastic_cfi}
			If $\stfilt_0^- \in \mathcal{X}$, (A,C) is observable and for all $t\geq n+1$ either \ref{as:larger_Q} with Assumption \ref{as:sufficient_feasibility} hold or \ref{as:stabilizing_condition} is satisfied, then there exists $s^\prime \geq 0$ such that
			\begin{equation}\label{e:converge_variance_cfi}
			\EE \left[\norm{\st_t - \stfilt_t^{\text{cf}}}^2 \right] \longrightarrow s^{\prime}.
			\end{equation}
		\end{theorem}	
	Proofs of theorems \ref{prop:FIE_stability}, \ref{prop:MHE_stability} and \ref{th:stochastic_cfi} are given in Appendix \ref{s:proofs_constrained}.
	\section{Numerical experiments}\label{s:experiments}
	For numerical experiments, we consider the benchmark model of a
	well-mixed, constant volume, isothermal batch reactor.  This
	model has previously been considered in~\cite{EKF_MHE, MHE2014_constrained}. The system dynamics is given by \eqref{e:process}, where
	\[
	A = \bmat{0.8831 & 0.0078 & 0.0022 \\ 0.1150 & 0.9563 & 0.0028
		\\ 0.1178 & 0.0102 & 0.9954}, \quad C = \bmat{32.84 & 32.84  & 32.84}
	\]
	The observability index of $(A, C)$ is $3$. The additive process and
	measurement noise are both assumed to be Gaussian with zero means, and
	variances, $(0.01)^2I$ and $(0.25)^2$, respectively.  
	The mean of the initial prior is $\stfilt_0^- =
	\bmat{1 &1 & 4}\transp$. Since the states represent concentration of
	chemicals in the batch reactor process, these cannot be negative.
	Therefore, the estimated states are constrained to lie in the set
	$\mathcal{X} \Let \{ x \in \R^d \mid x \geq 0 \}$. 
	
	\begin{example}\label{ex:Gaussian}
		In the first experiment, we assume that initial state is also Gaussian with prior mean $\stfilt_0^-$ and prior variance $\Sigma_0^- = I$. This is evident that simulated state of the system can be negative due to the presence of Gaussian noises in simulation but we consider this example for a fair comparison with minimum energy MHE (MEMHE) \cite{MHE01}. 
		\par We demonstrate a comparison between MEMHE and our proposed approach CMHE in Fig.\ \ref{fig:mse}. MEMHE is simulated by using \texttt{nmhe} object of freely available MATLAB based software package mpctools \cite{mpctools}, which is based on CasAdi \cite{casadi} and solver Ipopt \cite{Ipopt}. For CMHE, we use MATLAB-based software package YALMIP \cite{ref:lofberg-04} and a solver SDPT3-4.0 \cite{ref:toh-06} to solve the underlying optimization programs. We chose the optimization horizon $N=4$ for both approaches and simulated for $N_s = 1000$ sample paths. The empirical mean squared error $e_t$ for both approaches is computed by the following formula:
		\begin{equation}\label{e:mse}
		e_t = \frac{1}{N_s} \sum_{i=1}^{N_s}\norm{\st_t^i - \stfilt_t^i}^2 ,
		\end{equation}
		where $\st_t^i$ and $\stfilt_t^i$ denote the simulated and estimated states, respectively, at time $t$ in the $i^{\text{th}}$ path. 
		\par Fig.\ \ref{fig:mse} depicts that empirical mean squared error in our approach is smaller than that in MEMHE. Interestingly, at $t=0$ both approaches have almost same $e_t$ but in our approach it immediately drops by approximately one unit and keeps monotonically decreasing after then. However, in case of MEMHE a slight increase is observed at $t=2$ and after that it monotonically decreases but always remains higher than that of our approach.     
		
		\begin{figure}
		\begin{center}	
			\begin{adjustbox}{width = 0.8\columnwidth}
				\includegraphics{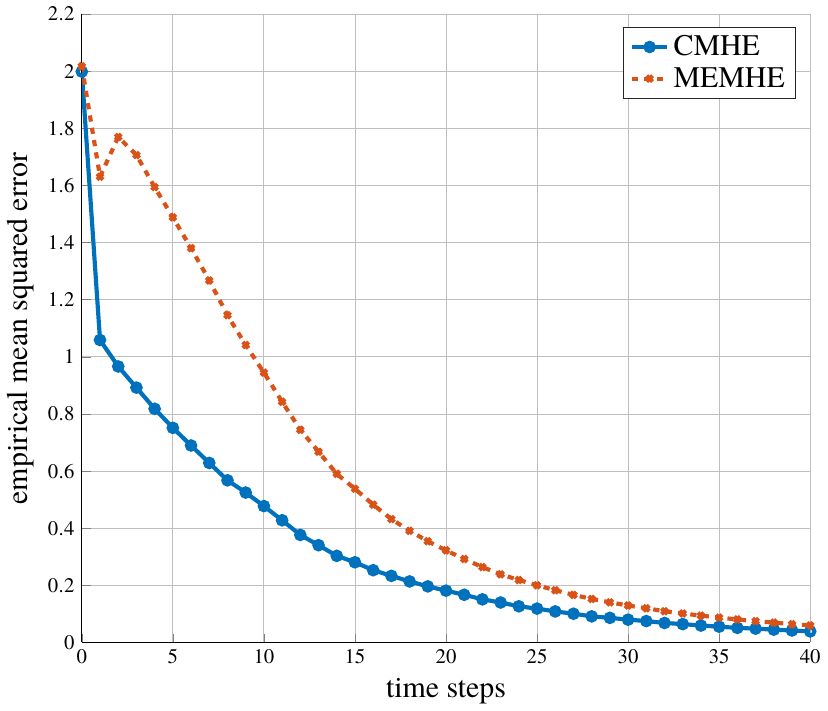}
			\end{adjustbox}	
		\end{center}
			\caption{The empirical mean squared error for 1000 sample paths is smaller in our proposed approach than that in standard MHE when initial state has Gaussian distribution.}
			\label{fig:mse}
		\end{figure}
	\end{example}
	\begin{example}\label{ex:uniform}
		In this experiment, we consider initial state to be uniformly distributed between $[0, 2\stfilt_0^-]$. Rest of the simulation data is same as in Experiment \ref{ex:Gaussian}. We simulate for $N_s = 100$ sample paths and compare between our proposed approach CMHE and standard MEMHE in Fig. \ref{fig:mse_uniform}. The empirical mean squared error $e_t$ is computed according to \eqref{e:mse}. Fig. \ref{fig:mse_uniform} depicts that both approaches have almost the same empirical mean squared error for 100 sample paths.  
		\begin{figure}
			\begin{center}
			\begin{adjustbox}{width = 0.8\columnwidth}
				\includegraphics{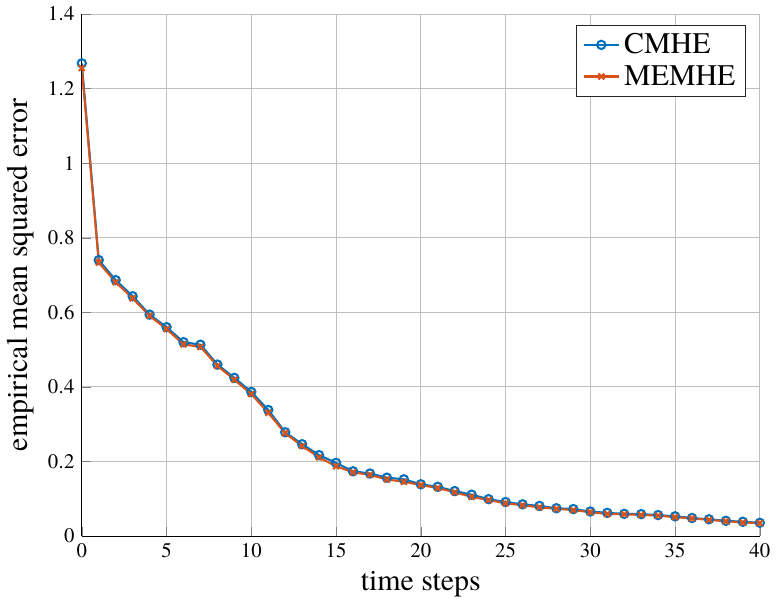}
			\end{adjustbox}	
		\end{center}
			\caption{The empirical mean squared error for 100 sample paths is almost same in our proposed approach and standard MHE when initial state has uniform distribution.}
			\label{fig:mse_uniform}
		\end{figure}
		
	\end{example}	
	\begin{example}
		In this experiment, we choose optimization horizon $N=3$ and simulate only for one sample path. Rest of the simulation data is same as in Experiment \ref{ex:uniform}. We compare the norm of estimate and cost by using CMHE and CFIE in Fig. \ref{fig:norm_cost}. Both approaches give almost same estimate and  incur almost same cost even though the optimization problem of CFIE has intermediate constraints, which are absent in CMHE.
		\begin{figure}
			\begin{center}
			\begin{adjustbox}{width = 0.8\columnwidth}
				\includegraphics{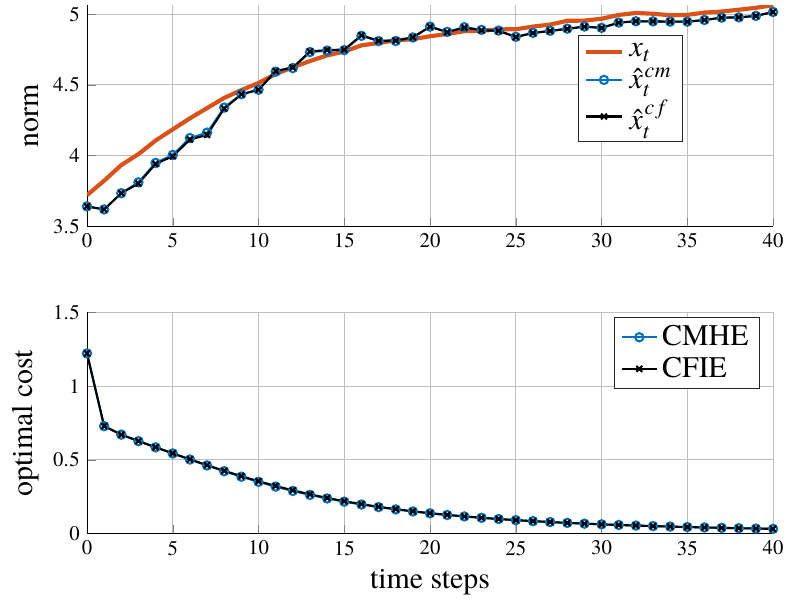}
			\end{adjustbox}	
		\end{center}
			\caption{Norm of estimate and optimal cost in both CMHE and CFIE are almost same.}
			\label{fig:norm_cost}
		\end{figure}
	\end{example}
	
	\section{Conclusions and directions for future research}\label{s:conclusion}
	In this paper, the minimum variance duality is used to convert the minimum
	variance estimation problem into a deterministic optimal control
	problem. The main contribution is the specification and the
	stability analysis of the FIE and MHE algorithms in the presence of
	state constraints. The proposed algorithms are distinct from and
	possess several useful features compared to the standard MHE
	algorithms based on the use of the minimum energy duality. In
	particular, there is no need to run a KF in parallel to approximate
	the terminal cost for the MHE. Both	the constrained FIE and MHE algorithms are stable in the sense of an	observer. Moreover, stochastic stability of constrained FIE
	is also established.	
	\par This work opens up several avenues for future research:
	Some ideas of linear model predictive control with time
	varying terminal cost and constraints ~\cite{timevarying_cost_constraint}, and approximate
	dynamic programming methods with accumulating constraints
	\cite{constrainedDP} may be useful for the further study
	of the constrained MHE. Several interesting extensions of the
	proposed approach may be possible including control design 
	\cite{Copp_Hespanha_simultaneous}, systems with intermittent
	observations \cite{PDQ_intermittent}, distributed architecture
	\cite{farina_distributedMHE}, the problem of unknown prior
	\cite{MHE02, MHE20} and inclusion of pre-estimating
	observer \cite{MHE03,MHE02, MHE18_metamorphic}. 

	\appendix

	\section{Proofs of \secref{s:unconstrained}}\label{s:proofs_unconstrained}
	
	\begin{proof}[Proof of Lemma \ref{lem:duality}]
		Since $z_0 - C\transp \alpha_0 = I$, we have 	
		\begin{align}
		& \st_t  = (z_0 - C\transp \alpha_0)\transp \st_t = z_0\transp \st_t - \alpha_0\transp C\st_t .\label{e:x_t}
		\end{align}
		By using the system dynamics \eqref{e:process} and the dual dynamics \eqref{e:dual_process}, we get
		\begin{equation}\label{e:canceling terms}
		\begin{aligned}
		z_i \transp \st_{t-i}  &= z_i\transp(A\st_{t-i-1} + \wnoise_{t-i-1}) \\
		z_{i+1}\transp \st_{t-i-1}  &= z_i\transp A\st_{t-i-1} + \alpha_{i+1}\transp C \st_{t-i-1}.
		\end{aligned}
		\end{equation}
		We substitute \eqref{e:canceling terms} in the expression of $z_0\transp \st_t$ as follows:  
		\begin{align*}	
		& z_0\transp \st_t  = \sum_{i=0}^{t-1} (z_i\transp \st_{t-i} - z_{i+1}\transp \st_{t-i-1}) + z_t\transp \st_0 \\	 
		& z_0\transp \st_t  = \sum_{i=0}^{t-1}\left( z_i\transp \wnoise_{t-i-1} - \alpha_{i+1}\transp C\st_{t-i-1} \right) + z_t\transp \st_0 .
		\end{align*}
		We further substitute $z_0\transp \st_t$ in \eqref{e:x_t} to get
		\begin{align*}
		&\st_t = \sum_{i=0}^{t-1}z_i\transp \wnoise_{t-i-1} - \sum_{i=0}^{t}\alpha_{i}\transp C\st_{t-i} + z_t\transp \st_0 \\		
		&\st_t = \sum_{i=0}^{t-1}z_i\transp \wnoise_{t-i-1} - \sum_{i=0}^{t}\alpha_{i}\transp(\meas_{t-i}-\mnoise_{t-i}) + z_t\transp \st_0.
		\end{align*}
		Further, we consider the estimate \eqref{e:estimator} and compute $ \EE [(\st_t - \stfilt_t)(\st_t - \stfilt_t)\transp]$ as follows:
		\begin{align*}	
		& \st_t - \stfilt_t = z_t\transp \left( \st_0 -  \stfilt_0^- \right) + \sum_{i=0}^{t-1}z_i\transp \wnoise_{t-i-1} + \sum_{i=0}^{t} \alpha_{i}\transp\mnoise_{t-i} \\
		& \EE [(\st_t - \stfilt_t)(\st_t - \stfilt_t)\transp] = z_t \transp \Sigma_{0}^- z_t + \sum_{i=0}^{t-1} z_i\transp Q z_i + \sum_{i=0}^{t}\alpha_i\transp R\alpha_i \\
		& \quad \quad = z_t \transp \Sigma_{0}^- z_t + \alpha_t\transp R\alpha_t + \sum_{i=0}^{t-1} \ell_i = \Gamma_0(z_t) + S_t(\alpha_{0:t+1}),	
		\end{align*}
		since the process noise, measurement noise and initial states are mutually independent.
		Therefore, $\EE \left[ \norm{\st_t - \stfilt_t}^2 \right] = \trace \left( \Gamma_0(z_t) + S_t(\alpha_{0:t+1}) \right)$, where $z_t$ is obtained by \eqref{e:dual_process} and $\stfilt_t$ is given by \eqref{e:estimator}.		
	\end{proof}	
	%%%%%%%%%%%%%%%%%%%%%%%%%%%%%%%%%%%%%%%%%%%%%%%%%%%%%%%%%%%%%%%%%%%%%%%%%%%%%%%%%%
	\begin{proof}[Proof of Lemma \ref{lem:uncsontrained_mhe}]
		At $t=0$, we compute 
		\begin{align*} 
		\Gamma_0(z_0) &+ S_0(\alpha_0) = z_0\transp \Sigma_0^-z_0 + \alpha_0\transp R\alpha_0 \\
		& = \Sigma_0^- + \alpha_0\transp (R + C\Sigma_0^- C\transp)\alpha_0 + \Sigma_0^-C\transp \alpha_0 + \alpha_0\transp C \Sigma_0^-. 	
		\end{align*}
		Since $\alpha_{0\mid 0} = \argmin \trace(\Gamma_0(z_0) + S_0(\alpha_0)) = -(C\Sigma_0^{-}C\transp + R)^{-1}C\Sigma_0^{-}$, due to our convention \eqref{e:variance_unconstrained} we obtain $\Sigma_0 = $ 
		\begin{equation}\label{e:first term}
		\Gamma_0(z_{0\mid 0})+S_0(\alpha_{0\mid 0}) = \Sigma_0^-   -  \Sigma_0^- C\transp (C\Sigma_0^- C\transp + R)^{-1}C\Sigma_0^-.
		\end{equation}  
		The FIE cost can be written as 
		\begin{align*}
		& \trace(\Gamma_0(z_t) + S_t(\alpha_{0:t+1})) = \trace\left(\Gamma_0(z_t) + \alpha_{t}\transp R \alpha_{t} + \sum_{i=0}^{t-1} \ell_i  \right) \\
		&= \trace\left( z_t\transp \Sigma_0^-z_t + \alpha_t\transp R\alpha_t + z_{t-1}\transp Qz_{t-1} + S_{t-1}(\alpha_{0:t}) \right).
		\end{align*}
		We substitute $z_t = A\transp z_{t-1} + C\transp \alpha_t$ in the above expression and the minimizer $\alpha_{t\mid t}\transp = -z_{t-1\mid t}\transp A \Sigma_0^-C\transp (C\Sigma_0^- C\transp +R)^{-1}$. Further, by substituting $\Sigma_0$ from \eqref{e:first term}, we get
		\begin{align}\label{e:cost_oneminus}
		& \trace(\Gamma_0(z_t) + S_t(\alpha_{0:t+1})) = \trace\left( z_{t-1}\transp (A\Sigma_0 A\transp + Q  )z_{t-1} + S_{t-1}(\alpha_{0:t})\right) \notag \\
		& = \trace\left( z_{t-1}\transp \Sigma_1^- z_{t-1} + S_{t-1}(\alpha_{0:t})\right),
		\end{align}
		where the last equality is due to our definition \eqref{e:prior}. Therefore, $\Gamma_1(\cdot)$ can be written as $\Gamma_1(z_{t-1}) = z_{t-1}\transp \Sigma_1^- z_{t-1}$. The above expression of cost \eqref{e:cost_oneminus} at time $t=1$ gives $\Sigma_1 = z_{0\mid 1}\transp \Sigma_1^- z_{0\mid 1} + S_0(\alpha_{0\mid 1})$, where $\alpha_{0\mid 1} = -(C\Sigma_1^{-}C\transp + R)^{-1}C\Sigma_1^{-}$.
		By repeating the above process $t-N$ times, we obtain
		\begin{align*}
		\trace(\Gamma_0(z_t) + S_t(\alpha_{0:t+1} ) =  \trace\left( z_{N}\transp \Sigma_{t-N}^- z_{N} + S_N(\alpha_{0:N+1})\right), 
		\end{align*}
		and therefore, we can define $\Gamma_{t-N}(z_N) = z_{N}\transp \Sigma_{t-N}^- z_{N}$.  
		Now for $t\geq N > 0$, we consider the expression of $\stfilt_t$:
		\begin{align*}
		& \stfilt_t = z_{t \mid t}\transp \stfilt_0^- - \alpha_{t \mid t} \transp \meas_0 - \sum_{i=0}^{t-1}\alpha_{i \mid t}\transp \meas_{t-i},\text{ where} \\
		& z_{t \mid t}\transp \stfilt_0^- - \alpha_{t \mid t} \transp \meas_0 = z_{t-1 \mid t}\transp A\stfilt_0^- - \alpha_{t \mid t}\transp (\meas_0 - C\stfilt_0^-).
		\end{align*}
		By substituting $\alpha_{t\mid t}$ in the above expression, we get	$z_{t \mid t}\transp \stfilt_0^- - \alpha_{t \mid t} \transp \meas_0 = z_{t-1\mid t}\transp \left(A\stfilt_0^- + A \Sigma_0^- C\transp (C\Sigma_0^- C\transp + R)^{-1}(\meas_0 - C\stfilt_0^-)  \right) = z_{t-1 \mid t}\transp A \stfilt_0$, which implies $\stfilt_t = z_{t-1 \mid t}\transp \stfilt_1^- - \sum_{i=0}^{t-1}\alpha_{i\mid t}\transp \meas_{t-i}$,	
		where $\stfilt_1^- = A\stfilt_0$. At $t=1$, we can compute $\stfilt_1^-$ from the above expression. By repeating the above process $t-N$ times we obtain the desired expression \eqref{e:estimator_mhe}. 
	\end{proof}	
	\section{Proofs of \secref{s:constrained}}\label{s:proofs_constrained}
	\begin{lemma}\label{lem:monotonicity_sign_definite}
		If \ref{as:larger_Q} holds then $\trace(\Sigma_t^{\text{cf}}) \geq \trace(\Sigma_{t-1}^{\text{cf}}) $ for all $t \geq n+1$. 	
	\end{lemma}
	\begin{proof}[Proof]
		Let us define $s_t^\ast \Let \trace(\Sigma_t^{\text{cf}})$ for notational simplicity. The optimal cost at time $t$ by substituting $S_t(\alpha_{0:t+1 \mid t}) = \alpha_{t\mid t}\transp R\alpha_{t \mid t} + \sum_{i=0}^{t-1} \ell_{i \mid t}$ in \eqref{e:cfie_mean_variance} is given by
		\begin{equation}\label{e:present cost}
		s_t^\ast = \trace\left( \Gamma_0(z_{t\mid t}) + \alpha_{t\mid t}\transp R \alpha_{t\mid t} + \sum_{i=0}^{t-1}\ell_{i\mid t} \right),
		\end{equation}
		where $\ell_{i\mid t} \Let z_{i\mid t}\transp Q z_{i \mid t} + \alpha_{i \mid t}\transp R\alpha_{i \mid t}$.
		We can observe for all $t \geq n+1$ that the constraints \eqref{e:intermediate constraint} at time $t-1$ are same at time $t$ for $j = 0, \ldots , t-n-1$. Therefore, $\alpha_{0:t \mid t}$, the first $t$ number of decision variables computed at time $t$, is a feasible control sequence at time $t-1$. Due to the optimality of $\alpha_{0:t\mid t-1}$ at time $t-1$, we get the following inequality:  	 
		\begin{equation*}\label{e:bound on earlier cost}
		\begin{aligned}
		s_{t-1}^\ast &\leq \trace\left( \Gamma_0(z_{t-1\mid t}) + \alpha_{t-1\mid t}\transp R \alpha_{t-1\mid t} + \sum_{i=0}^{t-2}\ell_{i\mid t} \right) \\
		& = \trace\left(z_{t-1\mid t}\transp (\Sigma_0^- - Q) z_{t-1\mid t} + \sum_{i=0}^{t-1}\ell_{i\mid t} \right) \\
		& = \trace\left(z_{t-1\mid t}\transp (\Sigma_0^- - Q) z_{t-1\mid t} -\Gamma_0(z_{t\mid t}) - \alpha_{t\mid t}\transp R \alpha_{t\mid t} \right) + s_t^\ast, 
		\end{aligned}
		\end{equation*}
		where the last equality is obtained by substituting \eqref{e:present cost}.
		Since $Q -\Sigma_{0}^- \succeq 0 $, for all $t \geq n+1$ we get
		\begin{equation*}
		s_t^\ast - s_{t-1}^\ast \geq \trace \left(\Gamma_0(z_{t\mid t}) + \alpha_{t\mid t}\transp R \alpha_{t\mid t} + z_{t-1\mid t}\transp (Q -\Sigma_{0}^-)z_{t-1\mid t}\right) \geq 0. 
		\end{equation*}		
	\end{proof}
	%%%%%%%%%%%%%%%%%%%%%%%%%%%%%%%%%%%%%%%%%%%%%%%%%%%%%%%%%%%%%%%%%%%%%%%%%%%%%%%
	\begin{lemma}\label{lem:monotonicity_Lyapunov_condition}
		If \ref{as:stabilizing_condition} holds, then $\trace(\Sigma_t^{\text{cf}}) \leq \trace(\Sigma_{t-1}^{\text{cf}})$ for all $t \geq n+1$.
	\end{lemma}
		\begin{proof}[Proof]
		We can observe that $\alpha_{0:t\mid t-1}$ satisfies \eqref{e:intermediate constraint} at time $t$ for $j=1, \ldots , t-n$. We assumed that $\alpha_t = K_t z_{t-1\mid t-1}$ satisfies  \eqref{e:intermediate constraint} for $j=0$. Therefore, the control sequence $\alpha_{0:t\mid t-1}$ along with $\alpha_t = K_t z_{t-1\mid t-1}$ is a feasible control sequence at time $t$. We compute $z_t$ by substituting $\alpha_{0:t\mid t-1}$ and $\alpha_t$ in \eqref{e:dual_process}, which gives us  
		$z_t = A\transp z_{t-1\mid t-1} + C\transp K_t z_{t-1\mid t-1} =\tilde{A}z_{t-1\mid t-1}$. Now we recall the expression of the optimal cost $s_t^\ast \Let \trace(\Sigma_t^{\text{cf}})$ from \eqref{e:present cost}.
		The optimaliity of $\alpha_{0:t+1\mid t}$ in the presence of stability criterion \eqref{e:monotonicity_Lyapunov_condition} gives us
		\begin{align*}
		&s_t^\ast \leq  \trace\left( \Gamma_0(\tilde{A}_t z_{t-1\mid t-1}  ) + z_{t-1\mid t-1} \transp K_t\transp RK_t z_{t-1\mid t-1} + \sum_{i=0}^{t-1}\ell_{i\mid t-1}\right) \\
		&= \trace\left( z_{t-1\mid t-1} \transp ( \tilde{A}_t\transp \Sigma_0^- \tilde{A}_t + K_t\transp RK_t + Q )z_{t-1\mid t-1}\right) \\
		& \quad + \trace\left( \alpha_{t-1\mid t-1}\transp R\alpha_{t-1\mid t-1} + \sum_{i=0}^{t-2}\ell_{i\mid t-1} \right) \\
& \leq \trace \left( \Gamma_0(z_{t-1 \mid t-1}) + \alpha_{t-1\mid t-1}\transp R\alpha_{t-1\mid t-1} + \sum_{i=0}^{t-2}\ell_{i\mid t-1} \right) = s_{t-1}^\ast.
		\end{align*}  
	\end{proof}
	%%%%%%%%%%%%%%%%%%%%%%%%%%%%%%%%%%%%%%%%%%%%%%%%%%%%%%%%%%%%%%%%%%%%%%%%%%%%%%%%%%%%%%%%%%
	\begin{lemma}\label{lem:feasible_policy}
		If the Assumption \ref{as:positive_invariance} holds and the matrix pair $(A,C)$ is observable then there exists $s > 0$ such that for the system \eqref{e:nominal_system},  
		\begin{align*}
		\trace(\Sigma_t^{\text{cf}}) &\leq s \text{ for all } t \text{ and} \\
		\trace(\Sigma_t^{\text{cm}}) &\leq s \text{ for all } N \geq n \text{ for all } t.
		\end{align*}
	\end{lemma}	
\begin{proof}[Proof]
		Let us consider the expression of $z_t$ at $t=n$ from \eqref{e:dual_state}. We can write it in compact form: $z_n =  A^{\top n}(I +  C\transp \alpha_0) + \reachab_n(A\transp, C\transp)\alpha_{1:n}$.
		If we substitute 
		\begin{equation}\label{e:fasible at n}
		\alpha_{1:n} = -\reachab_n(A\transp, C\transp)^{\dagger}(A^{ \top n})(I + C\transp \alpha_0)
		\end{equation}
		in the above expression for some $\alpha_0 \in \R^{q\times d}$, we get $z_n = \zeros$. Now we consider the estimator \eqref{e:stfilt without z} and the nominal system \eqref{e:nominal_system}. By substituting $\meas_i = CA^i\st_0$ and $\st_t = A^t \st_0$ for the system \eqref{e:nominal_system} in \eqref{e:stfilt without z}, we get
		\begin{align}\label{e:error without noise}
		\stfilt_t  &= A^t (\stfilt_0^- - \st_0) + A^t\st_0 + \sum_{i=0}^t \alpha_{i}\transp C A^{t-i} (\stfilt_0^- - \st_0) \notag \\
		&= \st_t + \left( A^t + \sum_{i=0}^t \alpha_{i }\transp C A^{t-i} \right)(\stfilt_0^- - \st_0) \notag \\
		& = \st_t + z_{t}\transp (\stfilt_0^- - \st_0),
		\end{align} 
		where the last equality is due to \eqref{e:dual_state}. If we substitute $\alpha_{1:n}$ from \eqref{e:fasible at n} in the above expression at $t=n$, we get $\stfilt_n = \st_n \in \bar{\mathcal{X}} \subseteq \mathcal{X}$ because $z_n = \zeros$ under \eqref{e:fasible at n}. Therefore, \eqref{e:fasible at n} is feasible for \eqref{e:constrained_fie} at $t=n$. 
		Let us define 
		\begin{equation} 
		s^0 \Let \trace(\Gamma_0(z_n)+ S_n(\alpha_{0:n+1})),
		\end{equation}
		where $z_n$ and $S_n$ are obtained by applying the given policy \eqref{e:fasible at n}. \\ 
		For all $t\geq n$, define $\beta_{0:n+1} = \alpha_{0:n+1} $ and $\beta_i = \zeros$ for $i > n$. Under the policy $\beta_{0:t+1}$, we have $z_t = \zeros$ and therefore $\stfilt_t = \st_t$ for all $t\geq n$; this policy is feasible. Since  $\trace(\Gamma_0(z_t) + S_t(\beta_{0:t+1})) = s$, optimality of $\alpha_{0:t+1\mid t}$ gives $\trace(\Sigma_t^{\text{cf}}) \leq s^0 \text{ for all } t \geq n$.		
		For each $t \leq n-1$, $\trace(\Sigma_t^{\text{cf}}) \leq \trace \left( A^t \Sigma_0^- A^{t \top} + \sum_{i=0}^{t-1} A^i Q A^{i \top} \right)$ is bounded, where the inequality holds due to optimality of $\trace(\Sigma_t^{\text{cf}})$ and feasibility of $\alpha_{0:t+1} = \zeros$. Defining $s \Let \max \{\trace(\Sigma_0^{\text{cf}}), \trace(\Sigma_1^{\text{cf}}), \ldots, \trace(\Sigma_{n-1}^{\text{cf}}), s^0\}$, we get the first part of the result. 
		Similarly, we can observe that $\beta_{0:N+1}$ is feasible for \eqref{e:constrained_MHE} for all $N\geq n$ and $\Sigma_t^{\text{cm}} = \Sigma_{t}^{\text{cf}}$ for $t \leq N$.
	\end{proof}
	
	%%%%%%%%%%%%%%%%%%%%%%%%%%%%%%%%%%%%%%%%%%%%%%%%%%%%%%%%%%%%%%%%%%%%%%%%%%%%%
	\begin{proof}[Proof of Theorem \ref{prop:FIE_stability}]
		For any $t \geq 0$, the optimal cost $\trace(\Sigma_t^{\text{cf}}) \leq s$ due to Lemma \ref{lem:feasible_policy}. Therefore, \cite[Lemma 6]{trace_bound} gives us the bound
		$ \lambda_{\min}(\Sigma_0^-)\trace(z_{t\mid t} z_{t \mid t}\transp) \leq \trace(z_{t\mid t} \transp \Sigma_0^- z_{t\mid t}) \leq \trace(\Sigma_t^{\text{cf}}) \leq s$,
		which further implies
		\begin{equation}\label{e:bound on z}
		\trace(z_{t\mid t} z_{t\mid t}\transp) \leq \frac{s}{\lambda_{\min}(\Sigma_0^-)}.
		\end{equation} 
		Set $\norm{\stfilt_0^- - \st_0} < \delta$ and consider $\norm{\stfilt_t^{\text{cf}} - \st_t}^2$. Since from \eqref{e:error without noise} $\stfilt_t^{\text{cf}} - \st_t = z_{t\mid t}\transp (\stfilt_0^- - \st_0)$, by using the bound \eqref{e:bound on z} we get
		\begin{align}
		&\norm{\stfilt_t^{\text{cf}} - \st_t}^2 = (\stfilt_0^- - \st_0) \transp z_{t\mid t} z_{t \mid t}\transp (\stfilt_0^- - \st_0) \leq \lambda_{\max} (z_{t\mid t}  z_{t \mid t} \transp) \norm{\stfilt_0^- - \st_0}^2 \notag \\ & \leq  \trace(z_{t \mid t}  z_{t \mid t} \transp) \norm{\stfilt_0^- - \st_0}^2 \leq \frac{s}{\lambda_{\min}(\Sigma_0^-)}\delta^2 \teL \varepsilon^2. \label{e:error bound fie}
		\end{align}
		Therefore, for a given $\varepsilon > 0$, we can choose $\delta = \sqrt{\frac{\lambda_{\min}(\Sigma_0^-)}{s}}\varepsilon$ which results in $\norm{\stfilt_t^{\text{cf}} - \st_t} \leq \varepsilon$ when $\norm{\st_0 - \stfilt_0^-} \leq \delta$ for all $t\geq 0$.
		In order to prove convergence of $\stfilt_t^{\text{cf}}$ to $\st_t$ for the system \eqref{e:nominal_system}, we first consider the case when $Q - \Sigma_0^- \succeq 0$. For all $t\geq n+1$, $\trace(\Sigma_t^{\text{cf}})$ is a monotonically increasing sequence due to Lemma \ref{lem:monotonicity_sign_definite} and it is bounded above due to Lemma \ref{lem:feasible_policy}. Therefore, it is convergent. From Lemma \ref{lem:monotonicity_Lyapunov_condition}, $\trace(\Sigma_t^{\text{cf}}) - \trace(\Sigma_{t-1}^{\text{cf}}) \rightarrow 0$, which implies $\trace(z_{t\mid t}z_{t\mid t}\transp) \rightarrow 0$ because $\Sigma_0^- \succ 0$. Then \eqref{e:error bound fie} immediately confirms that $\norm{\stfilt_t^{\text{cf}} - \st_t } \rightarrow 0$ as $t \rightarrow \infty$. 
		Now, we consider the second case when the stabilizing condition \eqref{e:monotonicity_Lyapunov_condition} of Lemma \ref{lem:monotonicity_Lyapunov_condition} is satisfied (\ref{as:stabilizing_condition} holds). In this case, $\trace(\Sigma_t^{\text{cf}})$ is a monotonically decreasing sequence which is bounded below. Similar to the first case, the convergence of $\trace(\Sigma_t^{\text{cf}})$ implies $\trace(z_{t\mid t}z_{t\mid t}\transp) \rightarrow 0$, which further implies $\norm{\stfilt_t^{\text{cf}} - \st_t} \rightarrow 0$.		
	\end{proof}
		\begin{proof}[Proof of Theorem \ref{prop:MHE_stability}]
			Let us consider the expression of $\Sigma_t^{\text{cm}}$ from \eqref{e:constrained_estimator}, $\Sigma_t^{\text{cm}} = z_{N\mid t}\transp \Sigma_{t-N}^{\text{cm}-}z_{N\mid t} + \alpha_{N\mid t}\transp R \alpha_{N\mid t} + \sum_{i=0}^{N-1}\ell_{i\mid t}$, 
			where $\ell_{i\mid t} = z_{i\mid t}\transp Q z_{i\mid t} + \alpha_{i\mid t}\transp R \alpha_{i\mid t} $. By substituting the expression of $\Sigma_{t-N}^{\text{cm}-}$ from \eqref{e:prior constrained}, we get $\Sigma_t^{\text{cm}} = z_{N\mid t}\transp \left( A\Sigma_{t-(N+1)}^{\text{cm}}A\transp +Q  \right)z_{N\mid t} + \alpha_{N\mid t}\transp R \alpha_{N\mid t} + \sum_{i=0}^{N-1}\ell_{i\mid t} = \sum_{i=0}^{N}\ell_{i\mid t} + z_{N\mid t}\transp A\Sigma_{t-(N+1)}^{\text{cm}}A\transp z_{N\mid t}$.
			Let us define $\gamma_{t,j} \Let A\transp z_{N\mid t-(j-1)(N+1)}\gamma_{t,j-1}$ with $\gamma_{t,0} = I$. Therefore,
			\begin{equation}\label{e:compact_MHE_cost}
			\Sigma_t^{\text{cm}} = \sum_{i=0}^{N}\ell_{i\mid t} + \gamma_{t,1}\transp \Sigma_{t-(N+1)}^{\text{cm}}\gamma_{t,1}. 
			\end{equation}
			For any $t=k(N+1) +r$, where $k \in \N$ and $r \in \{0, \ldots , N\}$, define $V_j = \sum_{i=0}^N \ell_{i \mid t - j(N+1)}$, by recursively solving \eqref{e:compact_MHE_cost} we get:
			\begin{align}\label{e:compact_MHE_cost_solve}
			&\Sigma_t^{\text{cm}} = \sum_{j=0}^{k-1} \gamma_{t,j}\transp V_j \gamma_{t,j} + \gamma_{t,k}\transp \Sigma_r^{\text{cm}}\gamma_{t,k} \\
			&= \sum_{j=0}^{k-1} \gamma_{t,j}\transp V_j \gamma_{t,j} + \gamma_{t,k}\transp \left( z_{r\mid r}\transp \Sigma_0^- z_{r\mid r} + \alpha_{r\mid r}\transp R \alpha_{r\mid r} + \sum_{i=0}^{r-1} \ell_{i\mid r}  \right) \gamma_{t,k}. \notag
			\end{align} 
			Since $s \geq \trace(\Sigma_t^{\text{cm}})$ for $t = k(N+1) + r$ due to Lemma \ref{lem:feasible_policy}, after ignoring some non-negative terms, we get
$s \geq \trace(\Sigma_t^{\text{cm}}) \geq \trace(\gamma_{t,k}\transp z_{r\mid r}\transp \Sigma_0^- z_{r\mid r} \gamma_{t,k}) \geq \trace( \gamma_{t,k} z_{r\mid r} \gamma_{t,k}\transp z_{r\mid r}\transp)\lambda_{\min}(\Sigma_{0}^{-})$. Therefore, 
			\begin{equation} \label{e:bound on last MHE term}
			\trace( \gamma_{t,k} z_{r\mid r} \gamma_{t,k}\transp z_{r\mid r}\transp) \leq \frac{s}{\lambda_{\min}(\Sigma_{0}^{-})}.
			\end{equation}	
			Now, we consider the expression of estimator for CMHE for the system \eqref{e:nominal_system} and substitute the expression of $\stfilt_{t-N}^{\text{cm}-}$ according to our definition \eqref{e:prior constrained}. For $t= k(N+1)+r$, similar to \eqref{e:error without noise}, we consider $\stfilt_{t}^{\text{cm}} -\st_t = z_{N\mid t}\transp (\stfilt_{t-N}^{\text{cm}-} - \st_{t-N})
			= z_{N\mid t}\transp A(\stfilt_{t-(N+1)}^{\text{cm}} - \st_{t-(N+1)})
			=  \gamma_{t,k}\transp (\stfilt_r^- - \st_r) =  \gamma_{t,k}\transp z_{r\mid r}\transp (\stfilt_0^- - \st_0)$. Therefore,
			\begin{equation}\label{e:cmhe_estimation_error_nominal}
			\begin{aligned}
			& \norm{\stfilt_{t}^{\text{cm}} -\st_t}^2 = \norm{\gamma_{t,k}\transp z_{r\mid r}\transp (\stfilt_0^- - \st_0)}^2 \leq \lambda_{\max}(\gamma_{t,k}\transp z_{r\mid r}\transp)\norm{\stfilt_0^- - \st_0}^2 \\
			& \leq \frac{s}{\lambda_{\min}(\Sigma_{0}^{-})}\norm{\stfilt_0^- - \st_0}^2 \teL \varepsilon^2,
			\end{aligned}  
			\end{equation}
			where the last inequality is due to \eqref{e:bound on last MHE term}.
			Therefore, for a given $\varepsilon > 0$, we can choose $\delta = \sqrt{\frac{\lambda_{\min}(\Sigma_0^-)}{s}}\varepsilon$ which results in $\norm{\stfilt_t^{\text{cm}} - \st_t} \leq \varepsilon$ when $\norm{\st_0 - \stfilt_0^-} \leq \delta$ for all $N\geq n$ and $t\geq N+1$. This completes the first part of the proof. For the second part, we consider \eqref{e:compact_MHE_cost_solve} and take limit $t \rightarrow \infty$, we get
			\begin{equation*}
			\lim_{t \rightarrow \infty} \trace( \Sigma_t^{\text{cm}} ) = \lim_{k \rightarrow \infty } \trace\left(  \sum_{j=0}^{k-1} \gamma_{t,j}\transp V_j \gamma_{t,j} + \gamma_{t,k}\transp \Sigma_r^{\text{cm}}\gamma_{t,k} \right) \leq s,
			\end{equation*}		 
			which results in $\trace(\gamma_{t,k-1}\transp V_{k-1} \gamma_{t,k-1}) \rightarrow 0$ as $k \rightarrow \infty$. By substituting $V_{k-1} = \sum_{i=0}^N \ell_{i \mid t - (k-1)(N+1)}$, we conclude that $\trace(\gamma_{t,k-1}\transp  \ell_{N \mid N+1+r} \gamma_{t,k-1}) \rightarrow 0$ and therefore, $\trace(\gamma_{t,k-1}\transp  z_{N \mid N+1+r}\transp Q z_{N \mid N+1+r} \gamma_{t,k-1}) \rightarrow 0$. Since $Q \succ 0$, we get $\trace(z_{N \mid N+1+r} \gamma_{t,k-1} \gamma_{t,k-1}\transp  z_{N \mid N+1+r}\transp) \rightarrow 0$ as $k \rightarrow \infty$. Now, we consider the expression \eqref{e:cmhe_estimation_error_nominal} and substitute $\gamma_{t,k} =A\transp z_{N\mid N+1+r}\gamma_{t,k-1}$ to get $\stfilt_{t}^{\text{cm}} -\st_t =  (A\transp z_{N\mid N+1+r}\gamma_{t,k-1})\transp z_{r\mid r}\transp (\stfilt_0^- - \st_0)$.
			Since $s \geq \trace(\Sigma_r^{\text{cf}}) \geq \trace(z_{r\mid r}\transp \Sigma_0^- z_{r\mid r})$, we get $\norm{z_{r\mid r}}_F \leq \sqrt{\frac{s}{\lambda_{\min}(\Sigma_0^-)}}$. We have
			\begin{align*}
			&\norm{\stfilt_{t}^{\text{cm}} -\st_t} = \norm{(z_{r\mid r} A\transp z_{N\mid N+1+r}\gamma_{t,k-1})\transp  (\stfilt_0^- - \st_0)} \\
			& \leq \norm{z_{N \mid N+1+r} \gamma_{t,k-1}}_F \norm{z_{r\mid r}}_F \norm{A}_F \norm{\stfilt_0^- - \st_0} \\
			& \leq \norm{z_{N \mid N+1+r} \gamma_{t,k-1}}_F \sqrt{\frac{s}{\lambda_{\min}(\Sigma_0^-)}} \norm{A}_F \norm{\stfilt_0^- - \st_0}, 
			\end{align*}
			which implies $\norm{\stfilt_{t}^{\text{cm}} -\st_t} \rightarrow 0$ because $\norm{z_{N \mid N+1+r} \gamma_{t,k-1}}_F \rightarrow 0$ as $t \rightarrow \infty$. This completes the second part of the proof. 
		\end{proof}	
		\begin{proof}[Proof of Theorem \ref{th:stochastic_cfi}]
			If $Q - \Sigma_0^- \succeq 0$, $\trace(\Sigma_t^{\text{cf}})$ is a monotonically increasing sequence due to Lemma \ref{lem:monotonicity_sign_definite}. We get a feasible control sequence due to Assumption \ref{as:sufficient_feasibility}. Therefore, due to optimality $\trace( \Sigma_t^{\text{cf}}) \leq \trace\left(\Gamma_0(z_t) + \alpha_t\transp R\alpha_t + \sum_{i=0}^{t-1} \ell_i  \right)$,
			where $\ell_i = z_i\transp Q z_i + \alpha_i\transp R \alpha_i, \alpha_{i+1} = K_iz_i$ and $z_{i+1} = (A\transp + C\transp K_i)z_i$. Due to the Assumption \ref{as:sufficient_feasibility}, we have $\ell_0 \leq \lambda_0 Q$, and for $t\geq 1$, $\ell_t \leq \lambda_t z_{t-1}\transp Q z_{t-1} \leq \lambda_t \ell_{t-1} \leq \lambda_t \lambda_{t-1} \ldots \lambda_0 Q$. Let us define $\rho_t \Let \lambda_t \lambda_{t-1} \ldots \lambda_0$, then
			\begin{align*}
			\trace( \Sigma_t^{\text{cf}}) & \leq \trace\left(\Gamma_0(z_t) + \ell_0 + \alpha_t\transp R\alpha_t + \sum_{i=1}^{t-1} \ell_i \right) \\
			& \leq \trace\left( z_t\transp (\Sigma_0^- - Q)z_t + \ell_0 + \sum_{i=1}^{t}\ell_i \right)  \\
			& \leq \trace\left(\lambda_0 Q +  \sum_{i=1}^{t} \ell_i \right) \leq \trace\left( \sum_{i=0}^{t} \rho_i Q \right)
			= \trace(Q)\sum_{i=0}^{t} \rho_i. 
			\end{align*}
			Since $\frac{\rho_{i+1}}{\rho_i}  = \lambda_{i+1} < 1$ for each $i$, there exists $\bar{\rho} > 0$ such that $\sum_{i=0}^{t} \rho_i < \bar{\rho}$ for each t. Therefore, $\trace( \Sigma_t^{\text{cf}}) \leq \bar{\rho}\trace(Q) $ for each $t$. 
			Since $\trace(\Sigma_t^{\text{cf}})$ is a monotionically increasing sequence and is bounded above, there exists some $s^\prime > 0$ such that \eqref{e:converge_variance_cfi} holds. This completes the first part of the proof. \\
			For the second case, the stabilizing condition of Lemma \ref{lem:monotonicity_Lyapunov_condition} is satisfied, and $\trace(\Sigma_t^{\text{cf}}) \geq 0$ is monotonically decreasing for all $t \geq n+1$. Therefore, there exists some $s^\prime \geq 0$ such that \eqref{e:converge_variance_cfi} holds. This completes the second part of the proof.
		\end{proof}						
	\section*{Acknowledgment}
	This work was supported in part by Navy N00014-19-1-2373 and NSF 1739874. The first author is thankful to Jin W. Kim for suggesting a useful reference.

	\bibliographystyle{IEEEtran}        
	\bibliography{estimator_controller} 
	
\end{document}